# Random walk on a polygon

**Jyotirmoy Sarkar**[1]

*Indiana University Purdue University Indianapolis*

**Abstract:** A particle moves among the vertices of an $(m+1)$-gon which are labeled clockwise as $0, 1, \ldots, m$. The particle starts at 0 and thereafter at each step it moves to the adjacent vertex, going clockwise with a known probability $p$, or counterclockwise with probability $1 - p$. The directions of successive movements are independent. What is the expected number of moves needed to visit all vertices? This and other related questions are answered using recursive relations.

## 1. Introduction

Consider a particle subjected to a random walk over the vertices of an $(m+1)$-gon which are labeled $0, 1, 2, \ldots, m$ in the clockwise direction. Initially the particle is at 0 and thereafter at each step the particle moves to one of the adjacent vertices, going in the clockwise direction with probability $p$, or in the counterclockwise direction with probability $q = 1 - p$. Throughout we assume that $0 < p < 1$ is a known constant. We also assume that the directions of successive movements are independent of one another. We answer the following questions in this paper:

1. What is the probability that all vertices have been visited before the particle returns to 0 for the first time?
2. What is the probability that the last vertex visited is $i$ $(i = 1, 2, \ldots, m)$?
3. What is the expected number of moves needed to visit all vertices?
4. What is the expected additional number of moves needed to return to 0 after visiting all vertices?

Question 1 appears in [6] (page 234, Exercise 46), while Questions 2 and 4 appear in [5] (page 224, Exercises 4.27 and 4.28). For the symmetric case of $p = 1/2$, Question 2 is solved in [6] (page 80) and Question 3 in Daudin [2], who furthermore gives the first five moments of the number of moves needed to visit all vertices. The asymmetric case of Question 3 also has been studied in the literature, but by using more advanced techniques (described in the next two paragraphs below), and the recorded expressions for the expected time to visit all vertices are formidable in appearance. We, on the other hand, present an elementary method of solution and derive a simple expression that is easy to comprehend. We hope that our solution will attract a wider readership.

Regarding the asymmetric case of Question 3, Vallois [8] studies the probability distribution of the time (called the *cover time*) taken by an asymmetric random walk on the integer lattice starting at 0 to visit $m + 1$ distinct integers (including the initial 0). That this problem indeed solves the problem of asymmetric random walk on a polygon is seen easily by simply joining the two extrema with a direct edge. Vallois uses the Martingale technique to obtain the probability generating function (p.g.f.) of the cover time as the ratio of two polynomials. Thereafter, he obtains the expected cover time by differentiation and the probability masses by inversion of the p.g.f.

A most comprehensive study on the joint distribution of the cover time, the last vertex visited (which is an extremum of the set of distinct numbers) and the time taken to visit the other extremum is presented by Chong, Cowan and Holst [1]. They follow the same approach as in Feller's classical treatise [3]. They express the joint p.g.f. in terms of trigonometric functions with complex arguments (as done by Feller) and also alternatively in terms of hyperbolic functions with real arguments. Thereafter, they obtain marginal distributions and moments of each of these component random variables.

In this paper we resolve Questions 1–3 by dissecting each question into parts that resemble the setup of the celebrated Gambler's Ruin Problem. We identify a clockwise movement of the particle with the gambler winning a dollar, and a counterclockwise movement with her losing a dollar. To solve these various parts we construct appropriate recursive relations using conditional probability and Bayes' rule, and then we solve the recursive relations either via difference equations or by mathematical induction. Finally, we assemble the solutions to the component parts to answer the entire original question. We solve the symmetric case of $p = 1/2$ first, because of its simplicity, and leave it to the reader to verify that this solution can be obtained also by taking the limit as $p \to 1/2$ in the corresponding result for the asymmetric case. We also present graphs to illustrate our results.

In Section 2 we review the Gambler's Ruin Problem by first summarizing the well known results, and then solving a new question in that setup. Section 3 contains the answers to Questions 1 through 3. Section 4 presents, without proof, the answer to Question 4 and also to the question of long run proportion of visits to each vertex using the well known limit theorem for finite Markov chains.

## 2. Gambler's ruin problem

A gambler plays a series of independent bets wagering one dollar on each bet. She either wins one dollar (and gets back her wager) with probability $p$, or loses her wager with probability $q = 1 - p$. Initially she has $i$ dollars. She must quit broke when her fortune reaches 0 (no credit is allowed). Also, she has adopted a policy to quit a winner when her fortune reaches $N$ dollars. Here, $N \geq i$ is a predetermined integer where remains fixed throughout.

The familiar questions that are presented in many standard textbooks are: (a) What is the probability that the gambler quits a winner with $N$ dollars rather than goes broke? (b) What is the expected number of bets she plays until the game ends? Below we present a synopsis of the solutions to these questions in the form of Propositions 2.1 and 2.2. For details the reader may see [3, 5, 6, 7], for example. Thereafter, we answer a new question in the same context: (c) Given that the gambler quits a winner with fortune $N$, what is the conditional expected number of games she has played? Likewise, given that the gambler goes broke, what is the conditional expected number of games she has played?

### 2.1. Probability of quitting a winner

**Proposition 2.1.** *A gambler starts with $i$ dollars, wagers one dollar per bet and wins a dollar with probability $p$ in each bet. The probability, $P_{i:N}$, that her fortune*

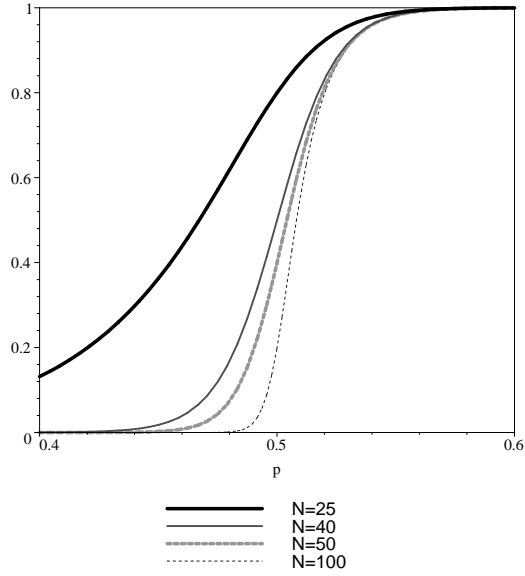

FIG 1. *The probability of reaching $N$ before $0$, starting from $i = 20$. Top to bottom the curves correspond to $N = 25, 40, 50, 100$.*

*reaches $N$ dollars before she goes broke is given by*

$$P_{i:N} = \begin{cases} \frac{i}{N} & \text{if } r = 1 \\ \frac{r^i - 1}{r^N - 1} & \text{if } r \neq 1, \end{cases} \tag{2.1}$$

*where $r = q/p = 1/p - 1$ is the odds of losing a bet.*

*Proof.* We condition on the outcome of the first bet. Let $F$ denote the event that the gambler wins the first bet and $W$ denote the event that the gambler quits a winner. Then we have the following recursive relations

$$\begin{aligned} P_{i:N} = P\{W\} &= P\{F\}\, P\{W|F\} + P\{F^c\}\, P\{W|F^c\} \\ &= p\, P_{i+1:N} + q\, P_{i-1:N} \quad \text{for } 1 \leq i \leq N-1, \end{aligned} \tag{2.2}$$

and the associated boundary conditions

$$P_{0:N} = 0 \quad \text{and} \quad P_{N:N} = 1, \tag{2.3}$$

which follow naturally from the quitting policy.

It is easy to verify that (2.1) satisfies (2.2) with the boundary conditions (2.3). The actual derivation of (2.1) is well known and may be found in any standard textbook such as [3, 5, 6, 7]. □

**Remark 1.** Note that $r = 1$ if and only if $p = 1/2$ (the game is fair). It should be pointed out that here and throughout the rest of the paper in all Propositions, Theorems and Corollaries the results are simpler in the case of $r = 1$, and they serve as a benchmark to verify the accuracy of the corresponding result in the asymmetric case of $r \neq 1$ by simply taking its limiting value as $r = q/p$ tends to 1. The continuity at $r = 1$ of expressions in (2.1), (2.4), (2.9), (2.15), (3.1), (3.3), (3.6), (3.7) and (4.1) can be verified by factoring $(r - 1)$ out and/or by applying L'Hospital's rule once or twice.

**Remark 2.** The gambler's expected fortune when the game ends, $N P_{i:N}$, decreases in $N \geq i$. Therefore, her optimum choice is $N = i$; that is, not to gamble at all. However, if she must gamble, she can minimize her loss by wagering the largest amount in each bet consistent with her quitting policy. See Feller [3].

### 2.2. Expected number of bets until the game ends

We only state the recursive relations needed in the elementary proof of Proposition 2.2 given below. Alternatively, for an elegant (though advanced) proof of Proposition 2.2 using Wald's Identities, see [5] (page 188), for example. For an elementary proof of Wald's Identities see [4].

**Proposition 2.2.** *In the setup of Proposition 2.1, the expected number of bets, $E_{i:N}$, until the gambler either reaches a fortune of $N$ dollars or goes broke, is given by*

$$(2.4) \qquad E_{i:N} = \begin{cases} i\,(N-i) & \text{if } r = 1 \\ \frac{r+1}{r-1}\left(i - N\frac{r^i-1}{r^N-1}\right) & \text{if } r \neq 1. \end{cases}$$

*Proof.* Again, by conditioning on the outcome of the first bet, we have

$$(2.5) \qquad E_{i:N} = 1 + p\,E_{i+1:N} + q\,E_{i-1:N} \quad \text{for } 1 \leq i \leq N-1,$$

and the quitting policy implies the boundary conditions

$$(2.6) \qquad E_{0:N} = 0 \quad \text{and} \quad E_{N:N} = 0.$$

It is easy to verify that (2.4) satisfies the system of equations (2.5)–(2.6). The details of the derivation are found in [6] (pages 234–235). □

**Remark 3.** The following two processes are equivalent: (1) Keep track of the gambler's fortune when the probabilities of success and failure in each bet are interchanged, and (2) Count by how many dollars the gambler's fortune is below $N$ when the final outcomes of quitting a winner and going broke are interchanged (but the probabilities remain unchanged). Therefore, we must have

$$(2.7) \qquad P_{i:N}(p \leftrightarrow q) = P_{i:N}(r \leftrightarrow r^{-1}) = 1 - P_{N-i:N}.$$

and

$$(2.8) \qquad E_{i:N}(p \leftrightarrow q) = E_{i:N}(r \leftrightarrow r^{-1}) = E_{N-i:N}.$$

Here and throughout the paper, the notation $f(p \leftrightarrow q)$ stands for a function identical in form $f(\cdot)$, but with arguments $p$ and $q$ interchanged. It is easily verified that (2.1) satisfies (2.7) and (2.4) satisfies (2.8).

### 2.3. Conditional expected number of bets, given that the gambler quits a winner

The question dealt with in this Subsection arises quite naturally in the context of the Gambler's Ruin problem. However, our literature search did not reveal its documentation anywhere; although we did find equation (2.13) below stated in [6] as Exercise 47 (page 234). The result of this Subsection is new, and so it is proved in complete detail.

**Theorem 2.1.** *In the setup of Proposition 2.1, the conditional expected number of bets, $W_{i:N}$, given that the gambler reaches a fortune of $N$ dollars before going broke, is given by*

$$(2.9) \qquad W_{i:N} = \begin{cases} \frac{1}{3}(N-i)(N+i) & \text{if } r = 1 \\ \frac{r+1}{r-1}\left[N\frac{r^N+1}{r^N-1} - i\frac{r^i+1}{r^i-1}\right] & \text{if } r \neq 1. \end{cases}$$

*Proof.* Note that for $i = 1 < N$, given that $W$ occurs, the gambler surely must have won the first bet, and thereafter she needs an expected number of $W_{2:N}$ additional bets to quit the game a winner. Hence, we have

$$(2.10) \qquad W_{1:N} = 1 + W_{2:N}.$$

Next, for $2 \leq i \leq N-1$, by conditioning on the outcome of the first bet, as in the proof of Proposition 2.2, we have

$$(2.11) \qquad W_{i:N} = 1 + P\{F|W\}\, W_{i+1:N} + (1 - P\{F|W\})\, W_{i-1:N}.$$

Lastly, from the quitting policy, we have the boundary condition

$$(2.12) \qquad W_{N:N} = 0.$$

Now, by Bayes' rule, for $1 \leq i \leq N-1$, we have

$$(2.13) \qquad P\{F|W\} = \frac{P\{F\}P\{W|F\}}{P\{W\}} = \frac{p\, P_{i+1:N}}{P_{i:N}} = \begin{cases} \frac{i+1}{2i} & \text{if } r = 1 \\ \frac{1}{r+1}\frac{r^{i+1}-1}{r^i-1} & \text{if } r \neq 1, \end{cases}$$

in view of Proposition 2.1.

Putting (2.13) in (2.11), one can verify that (2.9) satisfies the system of equations (2.10)–(2.12). The derivation of (2.9) is given in Appendix A. □

**Remark 4.** From (2.9) it is straightforward to verify that $W_{i:N}$ remains unchanged when the probabilities of success and failure in each bet are interchanged. That is,

$$(2.14) \qquad W_{i:N}(p \leftrightarrow q) = W_{i:N}(r \leftrightarrow r^{-1}) = W_{i:N}.$$

This is a pleasantly surprising result, especially in light of Remark 3.

**Corollary 2.1.** *In the setup of Proposition 2.1, the conditional expected number of bets, $B_{i:N}$, given that the gambler goes broke before reaching a fortune of $N$ dollars, is given by*

$$(2.15) \qquad B_{i:N} = \begin{cases} \frac{1}{3}i(2N-i) & \text{if } r = 1 \\ \frac{r+1}{r-1}\left[N\frac{r^N+1}{r^N-1} - (N-i)\frac{r^{N-i}+1}{r^{N-i}-1}\right] & \text{if } r \neq 1. \end{cases}$$

The Corollary follows from Theorem 2.1, since, in view of Remark 4,

$$(2.16) \qquad B_{i:N} = W_{N-i:N}(p \leftrightarrow q) = W_{N-i:N}(r \leftrightarrow r^{-1}) = W_{N-i:N}.$$

**Remark 5.** Since the game ends with the gambler either quitting a winner or going broke, we must necessarily have

$$E_{i:N} = P_{i:N}\, W_{i:N} + (1 - P_{i:N})\, B_{i:N},$$

which can be verified using (2.1), (2.9) and (2.15).

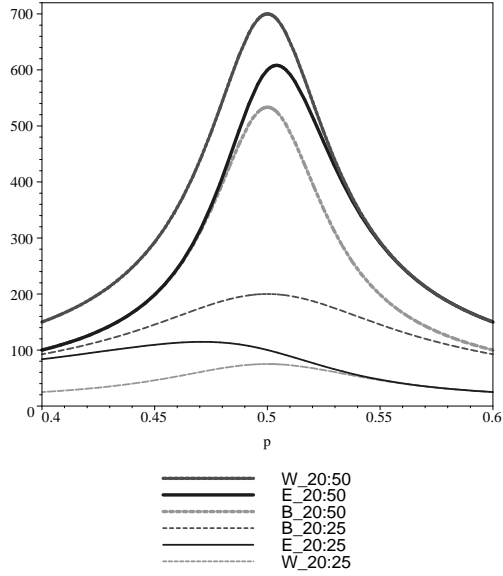

Fig 2. *The expected number of games until the gambler wins, or until the game ends, or until the gambler goes broke, starting from $i = 20$, for $N = 50, 25$. Top to bottom the curves correspond to $W_{20:50}, E_{20:50}, B_{20:50}, B_{20:25}, E_{20:25}, W_{20:25}$.*

## 3. Random walk on a polygon

In this Section we answer Questions 1–3 posed in Section 1. Throughout the paper we consider vertex $m + 1$ to be the same as vertex 0. In various steps of our solution we will renumber the vertices $m + 1 = 0, 1, 2, \ldots, m$, in that order, but the renumbering may be either in the clockwise or in the counterclockwise direction with an appropriate starting vertex. Also, we will dissect each question into parts that resemble the Gambler's Ruin Problem with the convention that a clockwise movement of the particle is identified with the gambler winning a dollar, while a counterclockwise movement is identified with her losing a dollar. Thus, for example, $F$ will denote the event that the first movement is clockwise (or the outcome of the first bet is a win), with $P\{F\} = p$.

### 3.1. Probability of visiting all vertices before returning to 0

**Theorem 3.1.** *Suppose that the vertices of an $(m + 1)$-gon are labeled $m + 1 = 0, 1, 2, \ldots, m$. A particle starting at 0 moves in each step to a neighboring vertex going clockwise with probability $p$, or counterclockwise with probability $q = 1 - p$. Let $r = q/p$. Let $A$ denote the event that the particle visits all vertices before its first return to 0. Then*

$$(3.1) \qquad P\{A\} = \begin{cases} \frac{1}{m} & \text{if } r = 1 \\ \frac{r-1}{r+1} \frac{r^m+1}{r^m-1} & \text{if } r \neq 1. \end{cases}$$

*Proof.* By conditioning on the first move, we have

$$(3.2) \qquad P\{A\} = p\, P\{A|F\} + q\, P\{A|F^c\}.$$

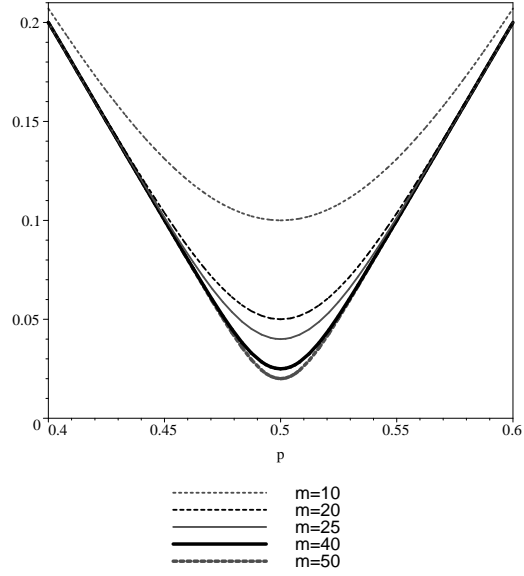

FIG 3. *The probability of visiting all m vertices before returning to* 0. *Top to bottom the curves correspond to* $m = 10, 20, 25, 40, 50$.

Note that $P\{A|F\}$ is the probability that the particle starting at vertex 1 visits vertex $m$ before vertex 0, and is equivalent to the probability that a gambler starting with 1 dollar quits a winner with $m$ dollars (before she goes broke). Hence, $P\{A|F\} = P_{1:m}$. Likewise, $P\{A|F^c\} = P_{1:m}(p \leftrightarrow q)$. Therefore, (3.2) becomes

$$P\{A\} = p\, P_{1:m} + q\, P_{1:m}(p \leftrightarrow q) = \begin{cases} \frac{1}{m} & \text{if } r = 1 \\ p\, \frac{r-1}{r^m - 1} + q\, \frac{r^{-1}-1}{r^{-m}-1} & \text{if } r \neq 1, \end{cases}$$

by Proposition 2.1, and simplifies to (3.1). □

**Remark 6.** In (3.1), note that $P\{A\} \geq 1/m$, with equality if and only if $r = 1$, implying that the probability of visiting all vertices before the particle returns to the origin is the smallest in the symmetric case of $p = 1/2$.

### 3.2. Probability distribution of the last vertex

**Theorem 3.2.** *In the setup of Theorem 3.1, let $L_i$ denote the event that the last vertex to be visited is vertex $i$. Then for $i = 1, 2, \ldots, m$,*

(3.3) $$P\{L_i\} = \begin{cases} \frac{1}{m} & \text{if } r = 1 \\ \frac{r^{m-i}(r-1)}{r^m - 1} & \text{if } r \neq 1. \end{cases}$$

*Proof.* Again, by conditioning on the first move, we have

(3.4) $$P\{L_i\} = P\{F\}\, P\{L_i|F\} + P\{F^c\}\, P\{L_i|F^c\} \\ = p\, P\{L_{i-1}\} + q\, P\{L_{i+1}\} \quad \text{for} \quad 2 \leq i \leq m-1,$$

and the boundary conditions

(3.5) $$P\{L_1\} = P_{1:m}(p \leftrightarrow q) \quad \text{and} \quad P\{L_m\} = P_{1:m},$$

by appropriate renumbering of vertices. To be more specific, by shifting all vertex labels one step in the counterclockwise direction, we have $P\{L_i|F^c\} = P\{L_{i+1}\}$ and $P\{L_m\} = P_{1:m}$. Likewise, by shifting vertex labels one step in the clockwise direction, we have $P\{L_i|F\} = P\{L_{i-1}\}$. Finally, by renumbering the vertices in the counterclockwise direction with the original vertex 1 becoming 0, we have $P\{L_1\} = P_{1:m}(p \leftrightarrow q)$.

In view of (3.5), we have the expressions for $P\{L_1\}$ and $P\{L_m\}$ from Proposition 2.1. It is straightforward to verify that (3.3) satisfies the system of equations (3.4)–(3.5). The derivation of (3.3) is similar to that of (2.1), and is left to the reader. □

**Remark 7.** We may rewrite (3.3) as

$$P\{L_i\} = \frac{r^{-i}}{r^{-1} + r^{-2} + \cdots + r^{-m}}.$$

Thus, the probability mass function for the last vertex is truncated geometric with support set $\{1, 2, \ldots, m\}$ and probabilities proportional to $r^{-i}$.

**Remark 8.** In the symmetric case of $r = 1$ all vertices are equally likely to be the last one visited. This result appears in [5] (page 80). It is indeed a mildly surprising result. When we surveyed the students in an undergraduate mathematics research class or in a graduate level stochastic processes class or even some mathematics/statistics faculty members, most of our subjects thought that the last vertex is more likely to be away from 0 than near 0.

### 3.3. Expected number of steps to visit all vertices

**Theorem 3.3.** *In the setup of Theorem 3.1, let $E[V]$ denote the expected number of steps needed to visit all vertices, and let $v_i = E[V|L_i]$ denote the conditional expected number of moves needed to visit all vertices given that the last vertex visited is $i$. Then for $i = 1, 2, \ldots, m$,*

(3.6) $$v_i = \begin{cases} \frac{1}{3}(m-1)(m+1) + (m+1-i)\,i & \text{if } r = 1 \\ v_i = \frac{r+1}{r-1}\left[m + i - 1 - \frac{2}{r-1} + \frac{2m}{r^m - 1} - (m+1)\frac{r^i - 1}{r^{m+1} - 1}\right] & \text{if } r \neq 1, \end{cases}$$

*and*

(3.7) $$E[V] = \begin{cases} \frac{m(m+1)}{2} & \text{if } r = 1 \\ \frac{r+1}{r-1}\left[m - \frac{1}{r-1} - \frac{m^2}{r^m - 1} + \frac{(m+1)^2}{r^{m+1} - 1}\right] & \text{if } r \neq 1. \end{cases}$$

*Proof.* If $L_1$ occurs, the particle first moves to vertex $m$, then it reaches vertex 2 before reaching vertex 1, and finally it reaches vertex 1. Similarly, if $L_m$ occurs, the particle starting at 0 first moves to vertex 1, then it reaches vertex $m-1$ before reaching vertex $m$, and finally it reaches vertex $m$. Hence, by appropriate renumbering of the vertices, we have

(3.8) $$v_1 = 1 + B_{m-2:m} + E_{1:m+1} \quad \text{and} \quad v_m = 1 + W_{2:m} + E_{m:m+1}.$$

Thus, $v_1$ and $v_m$ can be computed using (2.4), (2.9) and (2.15).

To compute $v_i$ for $2 \leq i \leq m-1$, we solve the following recursive relations

(3.9) $$v_i = 1 + q\, v_{i-1} + p\, v_{i+1} \quad \text{for } 2 \leq i \leq m-1.$$

To justify (3.9), note that by further conditioning on the outcome of the first bet, we have

$$v_i = E[V|L_i] = E[V|L_i \cap F]\, P\{F|L_i\} + E[V|L_i \cap F^c]\, P\{F^c|L_i\}$$
$$= v_{i-1}\, P\{F|L_i\} + v_{i+1}\, (1 - P\{F|L_i\})$$

by appropriate renumbering of the vertices. Also, by Bayes' rule, as in (2.13), we have

$$P\{F|L_i\} = \frac{P\{F\}\, P\{L_i|F\}}{P\{L_i\}} = p\, \frac{P\{L_{i-1}\}}{P\{L_i\}} = p\, r = q,$$

in view of Theorem 3.2.

One can verify that (3.6) satisfies (3.9) subject to (3.8). The derivation of (3.6) is given in Appendix B.

Finally, using $E[V] = \sum_{i=1}^m v_i\, P\{L_i\}$, (3.6) and Theorem 3.2, we get

$$E[V] = \begin{cases} \frac{1}{3}(m-1)(m+1) + \frac{1}{2}(m+1)^2 - \frac{1}{6}(m+1)(2m+1) & \text{if } r = 1 \\ \frac{r+1}{r-1}\left[m - \frac{1}{r-1} - \frac{m}{r^m - 1} + \frac{m+1}{r^{m+1}-1} - \frac{m(m+1)(r-1)r^m}{(r^m-1)(r^{m+1}-1)}\right] & \text{if } r \neq 1. \end{cases}$$

which simplifies to (3.7). □

**Remark 9.** When the probabilities of clockwise and counterclockwise movements of the particle are interchanged, we anticipate *a priori* the following relations to hold: (a) $P\{A\}(r \leftrightarrow r^{-1}) = P\{A\}$, (b) $P\{L_i\}(r \leftrightarrow r^{-1}) = P\{L_{m+1-i}\}$, (c) $v_i(r \leftrightarrow r^{-1}) = v_{m+1-i}$, and (d) $E[V](r \leftrightarrow r^{-1}) = E[V]$. The reader may verify these relations from (3.1), (3.3), (3.6) and (3.7).

## 4. Some further questions

In this Section, we answer Question 4 of Section 1, and a related question on the long run proportion of visits to each vertex. For the sake of brevity, the results are stated without proof.

### *4.1. Expected number of additional steps to return to the origin*

**Theorem 4.1.** *In the setup of Theorem 3.1, let $E[R]$ denote the expected number of additional steps needed to return to vertex* 0 *after visiting all vertices. Then, using* $E[R] = \sum_{i=1}^m P\{L_i\}\, E_{i:m+1}$, *we get*

$$E[R] = \begin{cases} \frac{1}{6}(m+1)(m+2) & \text{if } r = 1 \\ \frac{r+1}{r-1}\left[\frac{r}{r-1} - \frac{m(m+2)}{r^m-1} + \frac{(m+1)^2}{r^{m+1}-1}\right] & \text{if } r \neq 1. \end{cases}$$

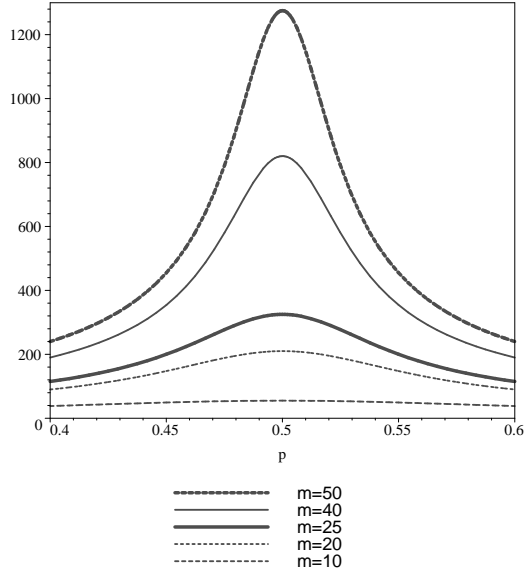

FIG 4. *The expected number of steps needed to visit all m vertices. Top to bottom the curves correspond to* $m = 50, 40, 25, 20, 10$.

### 4.2. Limiting results

When the random walk on the $(m+1)$-gon continues forever, we ask the following questions: (5) In the limit, how likely is the particle to be in each vertex? Or equivalently, what proportion of transitions enter into each vertex? (6) Starting from a particular vertex, what is the expected number of moves until the particle returns to that vertex?

To answer Questions 5 and 6, note that the polygonal random walk is a finite, irreducible and positive recurrent Markov chain. Hence, by invoking the well known limit theorem, (see, for example, [5], pages 175–177 and Problem 4.17 on page 221), we know that the limiting distribution is uniform; that is, the limiting probability for the particle to be in state $i$ is $1/(m+1)$, and the expected number of transitions needed to return to any state $i$ starting from itself is $m+1$, for $i = 0, 1, 2, \ldots, m$.

### Appendix A: Derivation of (2.9)

For the symmetric case of $r = 1$, substituting (2.13) in (2.11), we have

$$W_{i:N} = 1 + \frac{i+1}{2i} W_{i+1:N} + \frac{i-1}{2i} W_{i-1:N} \quad \text{for } 2 \le i \le N-1,$$

or equivalently,

$$(A.1) \qquad \frac{1}{2} i (i+1) (W_{i:N} - W_{i+1:N}) = i^2 + \frac{1}{2} (i-1) i (W_{i-1:N} - W_{i:N})$$

for $2 \le i \le N - 1$. This motivates us to define

$$(A.2) \qquad D_{i:N} = \frac{1}{2} i (i+1) (W_{i:N} - W_{i+1:N}) \quad \text{for } 1 \le i \le N - 1.$$

Then, (2.10) and (A.1) reduce to

$$D_{1:N} = 1 \quad \text{and} \quad D_{i:N} = D_{i-1:N} + i^2 \quad \text{for} \ 2 \leq i \leq N-1$$

which yields

(A.3) $$D_{i:N} = \sum_{j=1}^{i} j^2 = \frac{1}{6} i (i+1)(2i+1) \quad \text{for} \ 1 \leq i \leq N-1.$$

Next, substituting (A.3) in (A.2), we obtain

$$W_{i:N} - W_{i+1:N} = \frac{2 D_{i:N}}{i (i+1)} = \frac{2i+1}{3} \quad \text{for} \ 1 \leq i \leq N-1,$$

which, together with the boundary condition (2.12), yields

$$W_{i:N} = (W_{i:N} - W_{i+1:N}) + (W_{i+1:N} - W_{i+2:N}) + \cdots + (W_{N-1:N} - W_{N:N})$$
$$= \sum_{j=i}^{N-1} \frac{2j+1}{3} = \frac{1}{3}(N-i)(N+i).$$

This completes the derivation of (2.9) for the symmetric case.

For the asymmetric case of $r \neq 1$, substituting (2.13) in (2.11), we have

$$W_{i:N} = 1 + \frac{1}{r+1} \frac{r^{i+1}-1}{r^i-1} W_{i+1:N} + \frac{r}{r+1} \frac{r^{i-1}-1}{r^i-1} W_{i-1:N},$$

for $2 \leq i \leq N-1$, which we can rewrite as

$$\frac{1}{r+1} \frac{(r^{i+1}-1)}{(r^i-1)} (W_{i:N} - W_{i+1:N}) = 1 + \frac{r}{r+1} \frac{(r^{i-1}-1)}{(r^i-1)} (W_{i-1:N} - W_{i:N}),$$

or equivalently, as

(A.4)
$$\frac{(r^i-1)(r^{i+1}-1)}{(r+1)(r-1)^2} (W_{i:N} - W_{i+1:N})$$
$$= \left(\frac{r^i-1}{r-1}\right)^2 + r \frac{(r^{i-1}-1)(r^i-1)}{(r+1)(r-1)^2} (W_{i-1:N} - W_{i:N}),$$

for $2 \leq i \leq N-1$. Now letting

(A.5) $$C_{i:N} = \frac{(r^i-1)(r^{i+1}-1)}{(r-1)^2 (r+1)} (W_{i:N} - W_{i+1:N}) \quad \text{for} \ 1 \leq i \leq N-1,$$

we note that (2.10) and (A.4) reduce to

(A.6) $$C_{1:N} = 1 \quad \text{and} \quad C_{i:N} = r\, C_{i-1:N} + \left(1 + r + r^2 + \cdots + r^{i-1}\right)^2$$

for $2 \leq i \leq N-1$. By induction on $i \geq 1$, we solve (A.6) to get

(A.7) $$C_{i:N} = (r-1)^{-3} \left[r^{2i+1} - (2i+1)r^{i+1} + (2i+1)r^i - 1\right].$$

Substituting (A.7) in (A.5), after some algebraic simplification, we obtain for $1 \leq i \leq N-1$

$$W_{i:N} - W_{i+1:N} = \frac{r+1}{r-1} \left\{ 1 - \frac{2i}{r^i - 1} + \frac{2(i+1)}{r^{i+1} - 1} \right\},$$

which, together with the boundary condition (2.12), yields

$$\begin{aligned}
W_{i:N} &= (W_{i:N} - W_{i+1:N}) + (W_{i+1:N} - W_{i+2:N}) + \cdots + (W_{N-1:N} - W_{N:N}) \\
&= \frac{r+1}{r-1} \sum_{k=i}^{N-1} \left\{ 1 - \frac{2k}{r^k - 1} + \frac{2(k+1)}{r^{k+1} - 1} \right\} \\
&= \frac{r+1}{r-1} \left[ N - i - \frac{2i}{r^i - 1} + \frac{2N}{r^N - 1} \right],
\end{aligned}$$

completing the derivation of (2.9) for the asymmetric case.

**Appendix B: Derivation of (3.6)**

Letting $d_i = v_i - v_{i-1}$, we rewrite (3.9) as

(B.1) $$d_{i+1} = r \, d_i - (r+1), \quad \text{for} \quad i \geq 2$$

where $d_2$ is yet to be specified so that (3.8) holds.

For the symmetric case of $r = 1$, (B.1) simplifies to $d_i = d_2 - 2(i-2)$ for $i \geq 2$. Therefore, we have

(B.2) $$v_i = v_1 + d_2 + d_3 + \cdots + d_i = v_1 + (i-1)[d_2 - (i-2)].$$

By (3.8), (2.16) and (2.4), $v_m - v_1 = E_{m:m+1} - E_{1:m+1} = m - m = 0$. Hence, specializing (B.2) to $i = m$, we get $d_2 = m - 2$. Also by (3.8), (2.15) and (2.4), we have $v_1 = 1 + (m-2)(m+2)/3 + m$. Therefore, (B.2) reduces to

$$v_i = 1 + \frac{1}{3}(m-2)(m+2) + m + (i-1)(m-i) \quad \text{for} \quad 1 \leq i \leq m$$

which simplifies to (3.6).

For the asymmetric case of $r \neq 1$, (B.1) becomes

$$d_i = \left( d_2 - \frac{r+1}{r-1} \right) r^{i-2} + \frac{r+1}{r-1}.$$

Hence, we have

(B.3) $$v_i = v_1 + d_2 + d_3 + \cdots + d_i = v_1 + \left( d_2 - \frac{r+1}{r-1} \right) \frac{r^{i-1} - 1}{r - 1} + \frac{r+1}{r-1}(i-1).$$

Now, by (3.8), (2.16) and (2.4), we have

$$v_m - v_1 = E_{m:m+1} - E_{1:m+1} = \frac{r+1}{r-1}\left[(m-1) - (m+1)\frac{r^m - r}{r^{m+1} - 1}\right].$$

Therefore, specializing (B.3) to $i = m$, we get

(B.4) $$d_2 - \frac{r+1}{r-1} = \left[(v_m - v_1) - \frac{r+1}{r-1}(m-1)\right] \frac{r-1}{r^{m-1} - 1} = -\frac{(m+1)r(r+1)}{r^{m+1} - 1}.$$

Also by (3.8), (2.15) and (2.4), we have

$$\text{(B.5)} \qquad v_1 = \frac{r+1}{r-1}\left[m - \frac{2}{r-1} + \frac{2m}{r^m - 1} - (m+1)\frac{r-1}{r^{m+1} - 1}\right].$$

Substituting (B.5) and (B.4) in (B.3), and simplifying we establish (3.6) for the asymmetric case.

## Acknowledgments

The author wish to thank Benzion Boukai and Taiwo Salou for some discussions, Wenbo Li for some interesting conversations and David Aldous for pointing out [1]. The author is grateful to receive many thoughtful comments from two anonymous referees and an Associate Editor, which led to a significant improvement over an earlier draft. The research was partially funded by the Mathematical Association of America through its support of the IUPUI 2005 Undergraduate Research Experience in the Mathematical Sciences.